\documentclass[12pt]{amsart}
\usepackage{amsmath, amsthm, amscd, amsfonts}
\newtheorem{thm}{Theorem}[section]
 \newtheorem{cor}[thm]{Corollary}
 \newtheorem{lem}[thm]{Lemma}

 \numberwithin{equation}{section}

 \newfont{\Bbbb}{msbm10 scaled\magstephalf}
 \def\Fs{{F(p,q,s)}}
 \def\Nf{\| f\|_{F(p,q,s)}=|f(0)|+\{\sup_{a\in B}\int_B|\nabla f(z)|^p(1-|z|^2)^q g^s(z,a)dv(z)\}^{\frac{1}{p}}}
 \def\NBF{|f(0)|+\sup_{z\in B}(1-|z|^2)^\alpha|\nabla f(z)|}
 \def\nba{{\beta^\alpha}}
 \def\Xzu{\sup\limits_{\stackrel{u\in C^n-\{0\}}{\omega\in B_n}}|\psi(\omega)|\frac{(1-|w|^2)^\alpha}{(1-|\phi(w)|^2)^\frac{q+n+1}{p}}\Big\{\frac{(1-|\phi(w)|^2)|J_\phi(w)u|^2+|<\phi(w),J_\phi(w)u>|^2}{(1-|w|^2)|u|^2+|<w,u>|^2} \Big\}^{\frac{1}{2}}}
 \def\xzu{|\psi(\omega)|\frac{(1-|w|^2)^\alpha}{(1-|\phi(w)|^2)^\frac{q+n+1}{p}}\Big\{\frac{(1-|\phi(w)|^2)|J_\phi(w)u|^2+|<\phi(w),J_\phi(w)u>|^2}{(1-|w|^2)|u|^2+|<w,u>|^2} \Big\}^{\frac{1}{2}}}
 \def\Gp{G_{\frac{q+n+1}{p}}(\phi(\omega))|\nabla\psi(\omega)|(1-|\omega|^2)^\alpha}
 \def\wco{W_{\psi,\phi}}
 \def\nz{(1-|z|^2)}
 \def\na{(1-|a|^2)}
 \def\nw{(1-|\omega|^2)}
 \def\gs{g^s(z,a)}
 \def\nzw{|1-<z,w>|}
 \def\nza{|1-<z,a>|}
 \def\nua{|1-<u,a>|}
 \def\nu{(1-|u|^2)}
 \def\fwu{f_{\omega,u}}
 \def\fw{f_\omega}
 
 \def\gwu{g_{\omega,u}}
 \def\nfwus{\|f_{\omega,u}\|_{F(p,q,s)}}
 \def\ngwus{\|g_{\omega,u}\|_{F(p,q,s)}}
 \def\rw{r_\omega}
 \def\pw{\phi(\omega)}
 \def\npw{1-|\phi(\omega)|^2}
 \def\nrw{1-r_\omega^2}
 \def\Jpu{J_\phi(\omega)u}
 \def\nqp{\frac{n+1+q}{p}}
 \def\Ba{\Big{|}}
 \def\wu{\sqrt{\nw|u|^2+|<\omega,u>|^2}}
 \def\Uw{U_\omega}
 \def\nwa{\nw^\alpha}
 \def\npwq{(1-|\pw|^2)^\nqp}
 \def\Jpwu{\sqrt{(1-|\phi(w)|^2)|J_\phi(w)u|^2+|<\phi(w),J_\phi(w)u>|^2}}
 \def\be{\begin{equation}}
 \def\ee{\end{equation}}
 \def\bea{\begin{equation}\begin{array}{rcl}}
 \def\eea{\end{array}\end{equation}}

\begin{document}
\title[Weighted composition operators]{WEIGHTED COMPOSITION OPERATORS \\FROM $F(p,q,s)$
TO BLOCH TYPE SPACES\\
ON THE UNIT BALL}
\author[Z.H.Zhou and R.Y. Chen]{Zehua Zhou$^*$ \and  Renyu Chen}
\address{\newline Department of Mathematics\newline
Tianjin University
\newline Tianjin 300072\newline P.R. China.}

\email{zehuazhou2003@yahoo.com.cn}

\address{\newline Department of Mathematics\newline
Tianjin University
\newline Tianjin 300072\newline P.R. China.}
\email{fjcry@yahoo.com.cn}



\date{}
\thanks{\noindent $^*$ Zehua Zhou, Corresponding author. Supported in part by the National Natural Science Foundation of
China (Grand No.10371091), and LiuHui Center for Applied
Mathematics, Nankai University \& Tianjin University.}

\begin{abstract}
 Let $\phi(z)=(\phi_1(z),\cdots,\phi_n(z))$ be a
holomorphic self-map of $B$ and $\psi(z)$ a holomorphic function on
$B$, where $B$ is the unit ball of $\mbox{\Bbbb C}^n$. Let
$0<p,s<+\infty, -n-1<q<+\infty, q+s>-1$ and $\alpha\geq 0,$ this
paper gives some necessary and sufficient conditions for the
weighted composition operator $\wco$ induced by $\phi$ and $\psi$ to
be bounded or compact between the space $\Fs$ and $\alpha$-{\sl
Bloch} space $\beta^\alpha.$\end{abstract}

\maketitle

\section{Introduction}

 Let $dv$ be the $Lebesegue$ measure on the unit ball $B$ of
$\mbox{\Bbbb C} ^n$ normalized so that $v(B)=1$, and $d\sigma$ be
the normalized rotation invariant measure on the boundary $\partial
B$ of $B$ so that $\sigma(\partial B)=1.$
 $H(B)$ is the class of all holomorphic functions on $B$.

For $a\in B$, let $g(z,a)= \log|\varphi_a(z)|^{-1}$ be the {\sl
Green's} function on B with logarithmic singularity at $a$, where
$\varphi_a$ is the $M\ddot{o}bius$ transformation of $B$ with
$\varphi_a(0)=a, \varphi_a(a)=0, \varphi_a=\varphi_a^{-1}$.

Let $0<p,s<+\infty, -n-1<q<+\infty$ and $q+s>-1$. We say $f\in
F(p,q,s)$ provided that $f\in H(B)$ and
\begin{equation}
\Nf<+\infty \label{1}
\end{equation}
where
$$
\nabla f(z)=(\frac{\partial f(z)}{\partial z_1},\cdots,
\frac{\partial f(z)}{\partial z_n}).
$$

For the case $s=0,$ we say that $f$ belongs to the space
$F(p,q,0),$ if $$\|f\|_{F(p,q,0)}=|f(0)|+\left\{\int_B|\nabla
f(z)|^p(1-|z|^2)^q dv(z)\right\}^{\frac{1}{p}}<+\infty.$$

The corresponding "Little-oh" space $F_0(p,q,s)$ will be defined
as follows: $f$ belongs to $F_0(p,q,s)$, if
$$\lim\limits_{|a|\to 1}\int_B|\nabla f(z)|^p(1-|z|^2)^q
g^s(z,a)dv(z)=0.$$

For convenience, we also define $F_0(p,q,0)=F(p,q,0).$

For $\alpha \geq 0, f$ is said to be in the $Bloch$ space $\nba$
provided that $f\in H(B)$ and \be \|f\|_{\nba}=\NBF <+\infty.
\label{2} \ee As we all know, $\nba$ is a $Banach$ space when
$\alpha\geq 1$. The spaces $\beta^1$ and $\nba (0<\alpha<1)$ are
just the $Bloch$ space and the $Lipschitz$ spaces $L_{1-\alpha}$
respectively. From \cite{Yang/Ouyang} we know that a holomorphic
function $f\in\nba$ if and only if $\sup_{z\in
B}(1-|z|^2)^\alpha|Rf(z)|<+\infty,$ where \be Rf(z):=<\nabla
f(z),\bar{z}>=\Sigma^n_{j=1}z_j\frac{\partial f(z)}{\partial z_j}.
\label{3} \ee

Furthermore, by the Norm Equivalent Theorem we have \be
\|f\|_\nba\approx |f(0)|+\sup_{z\in B}(1-|z|^2)^\alpha|Rf(z)|,
\label{4} \ee where $M \approx N$ means the two quantities $M$ and
$N$ are comparable,
 that is there exist two positive constants $C_1$ and $C_2$ such that $C_1 M\leq N\leq C_2 M$.

If $p\geq 1,$ $\Fs$ is a $Banach$ space and it can be many
function spaces if we take some specific parameters of $p,q,$ and
$s$. For example, $\Fs=\beta^{\frac{q+n+1}{p}}$ when $s>n$;
$F(2,0,s)=Q_s, F(2,0,1)=BMOA$. It also includes Bergman spaces
$L^p_a=F(p,p,0)$ for $1\leq p<\infty$, Besove space
$B_p=F(p,p-n-1,0)$ for $n<p<+\infty$, Dirichlet space
$D_q=F(2,q,0)$ for $-1<q<\infty$, Hardy space $H^2=F(2,1,0).$ For
the definition of the spaces described above, we refer the reader
to see \cite{zh2, zhao}.

Let $\phi(z)=(\phi_1(z),\cdots,\phi_n(z))$ be a holomorphic self-map
of $B$ and $\psi(z)$ a holomorphic function on $B$, the composition
operator $C_\phi$ induced by $\phi$ is defined by
$$
(C_\phi f)(z)=f(\phi(z));
$$
the multiplication operator induced by $\psi$ if defined by
$$M_{\psi}f(z)=\psi(z)f(z);$$ and the weighted composition operator
$\wco$ induced by $\phi$ and $\psi$ is defined by
$$
(\wco f)(z)=\psi(z)f(\phi(z))
$$
for $z\in B$ and $f\in H(B)$. If let $\psi\equiv 1,$ then
$\wco=C_\phi$; if let $\phi=id$, then $W_{\psi,\phi}=M_{\psi}$. So
we can regard weighted composition operator as a generalization of a
multiplication operator and a composition operator. These operators
are linear.

Despite the simplicity of the definition of $W_{\phi,\psi}$ or
$C_{\phi}$, it is not uncommon for solutions of problems involving
this type of operator to require profound and interesting analytical
machinery; moreover, the study of composition operators has arguably
become a major driving force in the development of modern complex
analysis. The lucent texts \cite{Cowen/MacCluer}, \cite{Shapiro},
\cite{zh, zh2}, and conference proceedings collection \cite{cmr} are
good sources for information about many of the developments in the
theory of composition operators up to the middle of the last decade.

In the recent years, there have been many papers focused on studying
the composition operator in function spaces (say, for 1-dimensional
case, cf. \cite{Cowen/MacCluer} \cite{mm1} \cite{mm2}
\cite{Contreras/Hernandez-Diaz} \cite{Alfonso}), more recently,
[CSZ], \cite{Shi/Luo},\cite{Zhou/Shi1}, \cite{Zhou/Shi2},
\cite{Zhou1} discuss the composition operator in the Bloch space on
the unit polydisk. From those, they get some sufficient and
necessary conditions for composition operator to be bounded or
compact.

In this paper, our goal is to provide a theoretical characterization
of when $\phi$ and $\psi$ induce a  bounded or compact weighted
composition operators from $\Fs$ to $\alpha$-$Bloch$ space $\nba$.
Of course, if let $\psi=1$, it contains  the results about
composition operator; if let $\phi=id,$ it also contains the results
about multiplication operator. As some corollaries, we also give
some results between Bloch type space in the unit ball. These
results generalize previously known one-variable results obtained by
K. Madigan/A. Matheson \cite{mm2} and K. Madigan \cite{mm1},
respectively.

As far as the boundedness of  $\wco$ is concerned,  we can have the
following result:

\begin{thm}For $0<p,s<+\infty, -n-1<q<+\infty, q+s>-1,
\alpha\geq 0,$  $\phi$ be a holomorphic self-map of $B$, and $\psi$
a holomorphic function on $B$. Then the weighted composition
operator $\wco:\Fs\rightarrow\nba$ is bounded if and only if \bea
\Xzu <+\infty \label{5} \eea and \be \sup_{\omega\in B}\Gp<+\infty,
\label{6} \ee where $J_\phi$ is the Jabocian of $\phi$, and \bea
G_\nqp(\omega)=\left\{ \begin{array}{ll}
1, & 0<\nqp<1,\\
\log\frac{2}{1-|\omega|^2},  & \nqp=1, \\
\Big{(}\frac{1}{1-|\omega|^2}\Big{)}^{\nqp-1},& \nqp>1.
\end{array}\right.\label{7}
\eea
\end{thm}

On the other hand, the compactness of $\wco$ is much more
complicated than the case of boundedness. The result varies
sharply by the choosen of the real numbers $p$ and $q$.

\begin{thm}For $0<p,s<+\infty, -n-1<q<+\infty, q+s>-1,$  and $  \alpha\geq 0,$  $\phi$ be a holomorphic self-map of $B$,
$\psi$ a holomorphic function on $B$. If $ \frac{q+n+1}{p}\in
(0,1)$, then $\wco:\Fs\rightarrow\nba$ is compact if and only if
$\wco$ is bounded and \bea \sup\limits_{u\in\mbox{\Bbbb
C}^n-\{0\}}\xzu\rightarrow 0 \label{8} \eea when
$\phi(\omega)\rightarrow
\partial B$.
\end{thm}

\begin{thm}If $ \frac{q+n+1}{p}\geq 1,$ then
$\wco:\Fs\rightarrow\nba$ is compact if and only if $\wco$ is
bounded and \bea \sup\limits_{u\in\mbox{\Bbbb
C}^n-\{0\}}\xzu\rightarrow 0 \label{9} \eea \bea \Gp\rightarrow 0
\label{10} \eea when $\phi(\omega)\rightarrow \partial B$.
\end{thm}

Latter, in Lemma 2.1, Lemma 2.5, and Lemma 2.7, we will show that
$\Fs$ is just $\beta^{\nqp}$ when $s>n$, therefore by the theorems
above, we get

\begin{cor}Let $\phi=(\phi_1, \cdots, \phi_n)$ be a holomorphic self-map
of $B$ and $\psi$ a holomorphic function of $B, p>0, q>0$, then
$\wco:\beta^p\rightarrow\beta^q$ is bounded if and only if \bea
\sup\limits_{\stackrel{u\in \mbox{\Bbbb C}^n-\{0\}}{\omega\in
B_n}}|\psi(\omega)|\frac{(1-|w|^2)^q}{(1-|\phi(w)|^2)^p}\Big\{\frac{(1-|\phi(w)|^2)|J_\phi(w)u|^2+|<\phi(w),J_\phi(w)u>|^2}
{(1-|w|^2)|u|^2+|<w,u>|^2} \Big\}^{\frac{1}{2}}<+\infty \label{11}
\eea and \be \sup\limits_{\omega\in
B}G_p(\phi(\omega))|\nabla\psi(\omega)|(1-|\omega|^2)^q<+\infty .
\label{12} \ee
\end{cor}

\begin{cor}Let $\phi=(\phi_1, \cdots, \phi_n)$ be a
holomorphic self-map of $B$ and $\psi(z)$ a holomorphic function of
$B, 0<p<1, q>0$, then $\wco:\beta^p\rightarrow\beta^q$ is compact if
and only if $\wco$ is bounded and \bea \sup\limits_{u\in\mbox{\Bbbb
C}^n-\{0\}}|\psi(\omega)|\frac{(1-|w|^2)^q}{(1-|\phi(w)|^2)^p}\Big\{\frac{(1-|\phi(w)|^2)|J_\phi(w)u|^2+|<\phi(w),J_\phi(w)u>|^2}{(1-|w|^2)|u|^2+|<w,u>|^2}
\Big\}^{\frac{1}{2}}\rightarrow 0           \label{13} \eea when
$\phi(\omega)\rightarrow \partial B$.
\end{cor}

\begin{cor}Let $\phi=(\phi_1, \cdots, \phi_n)$ be a holomorphic self-map of $B$
and $\psi(z)$ a holomorphic function of $B, p\geq1, q>0$, then
$\wco:\beta^p\longrightarrow\beta^q$ is compact if and only if
 $\wco$ is bounded and
\bea \sup\limits_{u\in\mbox{\Bbbb
C}^n-\{0\}}|\psi(z)|\frac{(1-|w|^2)^q}{(1-|\phi(w)|^2)^p}\Big\{\frac{(1-|\phi(w)|^2)|J_\phi(w)u|^2+|<\phi(w),J_\phi(w)u>|^2}{(1-|w|^2)|u|^2+|<w,u>|^2}
\Big\}^{\frac{1}{2}}\rightarrow 0  \label{14} \eea \bea
G_p(\phi(\omega))|\nabla\psi(\omega)|(1-|\omega|^2)^q\rightarrow 0
\label{15} \eea when $\phi(\omega)\rightarrow \partial B$.
\end{cor}

The organization of this paper is as follows: We give some Lemmas in
Section 2, and prove Theorem 1.1, Theorem 1.2 and Theorem 1.3 of
this paper in Section 3, 4 and 5, respectively.

\section{Some Lemmas}

In the following, we will use the symbol $c$ or $C$ to denote a
finite positive number which does not depend on variables $
z,a,\omega$ and may depend on some norms and parameters $p, q, n,
\alpha, x, f$ etc, not necessarily the same at each occurrence.

In order to prove the main result, we will give some Lemmas first.

\begin{lem}If $0<p,s<+\infty, -n-1<q<+\infty, q+s>-1$, then
$\Fs\subset\beta^\nqp$ and $\exists c>0$ s.t. for $\forall f\in\Fs,
$ $\|f\|_{\beta^\nqp}\leq c\|f\|_\Fs. $
\end{lem}

\begin{proof} This lemma has been given by \cite{Zhang}, but for the convenience
of the reader, we will still give the proof here.

Suppose $f \in \Fs.$ Fixed $0<r_0<1,$ since $(Rf) \circ \varphi_a\in
H(B),$ so $|(Rf)\circ\varphi_a|^p$ is subharmonic in $B$. That is
\bea
|Rf(a)|^p=&|(Rf)\circ\varphi_a(0)|^p\leq\frac{1}{r_0^{2n}}\int_{r_0B}|(Rf)\circ\varphi_a(\omega)|^pdv(\omega)\\
          &=\frac{1}{r_0^{2n}}\int_{r_0B}|(Rf)|^p\frac{\na^{n+1}}{\nza^{(2n+2)}}dv(z).
\label{2.1}\eea From (5) in \cite{Zhuo/Ouyang}, we have
$$
\frac{1-r_0}{1+r_0}\na\leq\nz\leq\frac{1+r_0}{1-r_0}\nz
$$
as $z\in\varphi_a(r_0 B)$. Thus \be
\frac{\na^{n+1}}{\nza^{2n+2}\nz^qg^s(z,a)}\leq\frac{4^{n+1}}{\na^{n+1+q}}(\frac{1+r_0}{1-r_0})^{|q|}\log^{-s}\frac{1}{r_0}.
\label{2.2}\ee

 From (\ref{2.1}) and (\ref{2.2}), we get
$$
\begin{array}{ll}
&|Rf(a)|^p\leq\frac{1}{r_0^{2n}}\int_{\varphi_a(r_0B)}|Rf(z)|^p\frac{\na^{n+1}}{\nza^{2n+2}}dv(z)\\
&=\frac{1}{r_0^{2n}}\int_{\varphi_a(r_0B)}|Rf(z)|^p\nz^q\gs\frac{\na^{n+1}}{\nza^{2n+2}\nz^q\gs}dv(z)\\
&\leq\frac{4^{n+1}r_0^{-2n}}{\na^{n+1+q}}(\frac{1+r_0}{1-r_0})^{|q|}\log^{-s}\frac{1}{r_0}\|f\|^p_\Fs.
\end{array}
$$
This shows that $f\in\beta^{\frac{n+1+q}{p}}$ and
$\|f\|_{\beta^\nqp}\leq c\|f\|_\Fs. $
\end{proof}

\begin{lem}Let $p>0$, then there is a constant $c>0, $
for $\forall f\in\beta^p$ and $\forall z \in B,$ the estimate \be
|f(z)|\leq c G_p(z)\|f\|_{\beta^p}, \label{16} \ee holds, where the
function $G_\alpha$ has been defined at (\ref{7}).
\end{lem}

\begin{proof}For  $\forall f\in\beta^p(B_n),$ since
$||f||_{\beta^p}=|f(0)|+\sup\limits_{z\in B_n}(1-|z|^2)^p|\nabla
f(z)|$, we have
$$
|f(0)| \leq ||f||_{\beta^p},\quad \mbox{and}\quad |\nabla f(z)| \leq
\frac{||f||_{\beta^p}}{(1-|z|^2)^p}.
$$
but
$$
f(z)=f(0)+ \int_0^1<z,\overline{\nabla f(tz)}>dt.
$$
therefore
\begin{eqnarray*}
&&|f(z)|\leq |f(0)|+\int_0^1|z|\,|\nabla f(tz)|dt \\[6pt]
&\leq& ||f||_{\beta^p}+||f||_{\beta^p}\int_0^1
\frac{1}{(1-|tz|^2)^p}dt \leq
||f||_{\beta^p}\Big(1+\int_0^{|z|}\frac{dt}{(1-t^2)^p}\Big).
\end{eqnarray*}
when $p=1$,
$\displaystyle\int_0^{|z|}\displaystyle\frac{dt}{1-t^2}=\displaystyle\frac{1}{2}
\ln\displaystyle\frac{1+|z|}{1-|z|}\leq
\frac{1}{2}\displaystyle\ln\frac{4}{1-|z|^2}$, therefore
$$
|f(z)|\leq
\Big(1+\frac{1}{2}\ln\frac{4}{1-|z|^2}\Big)||f||_{\beta^p}.
$$

If $p\neq 1$, then
$$
\displaystyle\int_0^{|z|}\displaystyle\frac{dt}{(1-t^2)^p}=\displaystyle\int_0^{|z|}
\displaystyle\frac{dt}{(1-t)^p(1+t)^p}\leq
\displaystyle\int_0^{|z|}\displaystyle\frac{dt}{(1-t)^p}=\displaystyle\frac{1-(1-|z|)^{1-p}}{1-p},
$$
therefore when $0<p<1$, notice that
$\displaystyle\int_0^{|z|}\displaystyle\frac{dt}{(1-t^2)^p} \leq
\displaystyle\frac{1}{1-p} $ we get
$$
|f(z)|\leq \Big(1+\frac{1}{1-p}\Big)||f||_{\beta^p}.
$$
and when $p>1$
\begin{eqnarray*}
&&\int_0^{|z|}\frac{dt}{(1-t^2)^p}\leq
\frac{1-(1-|z|)^{1-p}}{1-p}\\
&&=\frac{1-(1-|z|)^{p-1}} {(p-1)(1-|z|)^{p-1}}\leq
\frac{2^{p-1}}{(p-1)(1-|z|^2)^{p-1}}
\end{eqnarray*}
so
$$
|f(z)|\leq
\Big(1+\frac{2^{p-1}}{(p-1)(1-|z|^2)^{p-1}}\Big)||f||_{\beta^p}.
$$
\end{proof}

\begin{lem}Let $\phi$ be a holomorphic self-map of $B$
and $\psi$ a holomorphic function on $B$, $K$ is an arbitrary point
set.
 Then $\wco :\Fs\rightarrow \nba$ is compact if and only if for any uniformly bounded sequence $\{f_{j,u}\}(j\in\mbox{\Bbbb N},u\in K)$
 in $\Fs$ which converges to zero uniformly for $u\in K$ and $z$ on any compact subsets of $B$ when  $j\rightarrow \infty$, $\|\wco f_{j,u}\|_\alpha\rightarrow 0$ holds.
\end{lem}

\begin{proof}Suppose that $\wco$ is compact and $\{f_{j,u}\}$ is a uniformly
bounded sequence in $\Fs$ which converges to zero uniformly for
$u\in K$ and $z$ on compact subsets of $B$ when $j\rightarrow
\infty$. Then $\{\wco(f_{j,u})\}$ has a subsequence
$\{\wco(f_{j_m,u_k})\}$ converges to $g\in\nba$. And by lemma 2.2,
$\forall$ compact subset $M\subset B,$ there is a positive constant
$C_M$ independent of $f_{j,u}$ such that
$$
|\psi(z)f_{j_m,u_k}(\phi(z))-g(z)|\leq C_M\|\psi f_{j_m,u_k}\circ
\phi-g\|_\alpha,
$$
for all $z\in M.$ Therefore $\{\psi(z)f_{j_m,u_k}(\phi(z))-g(z)\}$
converges to 0 uniformly on M. Notice that there is a constant
$c>0$, such that $|\psi(z)|<c, \forall z\in M$, and that $\phi(M)$
is compact in B, we have $$|\psi(z)f_{j_m,u_k}(\phi(z))|\leq
c|f_{j_m,u_k}(\phi(z))|\rightarrow 0, $$ uniformly on $M$. And by
the arbitrariness of $M$, we have $g\equiv 0$ on $B$. Since it is
true for arbitrary subsequence of $\{f_{j,u}\}$, we see that $\wco
f_{j,u}\to 0$ in $\nba$, when $j\to +\infty.$

On the other hand, fixed a point $b$ in $K$. Suppose
$\{g_{j,b}\}\subset K_r=B_\Fs(0,r)$, where $B_\Fs(0,r)$ is a ball in
$\Fs,$ then by lemma 2.2, $\{g_{j,b}\}$ is uniformly bounded in
arbitrary compact subset M of $B$. By Montel's Lemma, $\{g_{j,b}\}$
is regular, therefore there is a subsequence $\{g_{j_m,b}\}$ which
converges uniformly to $g_b\in H(B)$ on compact subsets of $B$. It
follows that $\nabla g_{j_m,b}\to\nabla g_b$ uniformly on compact
subsets of B.

Denote $B_k=B(0,1-\frac{1}{k})\subset \mbox{\Bbbb C}^n$, then
$$\begin{array}{ll}
&\int_B|\nabla g_b|^p\nz^q\gs dv(z)\\
&=\lim\limits_{k\to+\infty}\int_{B_k}\lim\limits_{m\to+\infty}|\nabla
g_{j_m,b}|^p\nz^q\gs dv(z)\\
&=\lim\limits_{k\to+\infty}\lim\limits_{m\to+\infty}\int_{B_k}|\nabla
g_{j_m,b}|^p\nz^q\gs dv(z).
\end{array}
$$
But $\{g_{j_m,b}\}\subset B_\Fs(0,r)$, then $$\int_{B_k}|\nabla
g_{j_m,b}|^p\nz^q\gs dv(z)<r^p,$$therefore
$$
\int_B|\nabla g_b|^p\nz^q\gs dv(z)\leq r^p.
$$
So $\|g_b\|_\Fs\leq r,$ and $g_b\in\Fs.$ Hence the sequence
$\{g_{j_m,b}-g_b\}$ is such that $\|g_{j_m,b}-g_b\|\leq 2r<\infty$
and converges to $0$ on compact
subsets of B, by the hypothesis of this lemma, we have that $$\psi
g_{j_m,b}\circ \phi\to\psi g_b\circ \phi$$ in $\nba$. Thus the set
$\wco(K_r)$ is relatively compact, finishing the proof.
\end{proof}

\begin{lem}Let $0< p< \infty,f\in H(B_n),$ then $f\in \beta^p(B_n)$ if and
only if
$$
\sup_{\stackrel{u\in \mbox{\Bbbb C}^n-\{0\}}{z\in
B_n}}\frac{(1-|z|^2)^p|\nabla f(z)u|}
{\sqrt{(1-|z|^2)|u|^2+|<z,u>|^2}} <+\infty.
$$
Furthermore
$$
||f||_{\beta^p}\approx |f(0)|+\sup_{\stackrel{u\in \mbox{\Bbbb
C}^n-\{0\}}{z\in B_n}}\frac{(1-|z|^2)^p|\nabla f(z)u|}
{\sqrt{(1-|z|^2)|u|^2+|<z,u>|^2}} .
$$
\end{lem}

\begin{proof}By Bergman metric in \cite{Timoney1}, there are constants $A_1>0$ and
$A_2>0$, such that
$$
\begin{array}{ll}
&\sup\limits_{u\in\mbox{\Bbbb C}^n-\{0\}}\frac{A_1|\nabla f(z)u|}{\sqrt{\nz|u|^2+|<z,u>|^2}}\leq |\nabla f(z)|\\
&\leq \sup\limits_{u\in\mbox{\Bbbb C}^n-\{0\}}\frac{A_2|\nabla
f(z)u|}{\sqrt{\nz|u|^2+|<z,u>|^2}}.
\end{array}
$$
Let $C_1=\min \{A_1,1\}$, $C_2=\min \{A_2,1\}$, we have  $$
\begin{array}{ll}
&C_1\{|f(0)|+\sup\limits_{u\in\mbox{\Bbbb C}^n-\{0\}}\frac{\nz^p|\nabla f(z)u|}{\sqrt{\nz|u|^2+|<z,u>|^2}}\}\leq\|f\|_{\beta^p}\\
&C_2\{|f(0)|+\sup\limits_{u\in\mbox{\Bbbb
C}^n-\{0\}}\frac{\nz^p|\nabla
f(z)u|}{\sqrt{\nz|u|^2+|<z,u>|^2}}\}.
\end{array}
$$
 This proof is completed.
\end{proof}

\begin{lem}For $0<p,s<+\infty, -n-1<q<+\infty, q+s>-1$, there exists $c>0$
such that \bea \sup\limits_{a\in
B}\int_B\frac{\nw^p}{\nzw^{n+1+q+p}}\nz^q\gs dv(z)\leq c, \label{18}
\eea
 for every $\omega \in B$ .
\end{lem}
\begin{proof}It's  easy to verify that if $0<x<\frac{1}{2}$, then
$-\log(1-x)<4x$. If we let $x=1-|\varphi_a(z)|^2$, then \bea
 &-\log(1-(1-|\varphi_a(z)|^2))<4(1-|\varphi_a(z)|^2),\eea or$$g(z,a)<2(1-|\varphi_a(z)|^2)
 $$ when $|\varphi_a(z)|^2>\frac{1}{2}$.
Note that $1-|\varphi_a(z)|^2=\frac{\na\nz}{\nza^2}$, therefore when
$s>n$ \be
\begin{array}{ll}
&\int_{\frac{1}{2}<|\varphi_a(z)|^2<1}\frac{\nw^p}{\nzw^{n+1+q+p}}\nz^q\gs dv(z)\\
&\leq\int_{\frac{1}{2}<|\varphi_a(z)|^2<1}\frac{\nw^p\nz^q 2^s\na^s\nz^s}{\nzw^{n+1+q+p}\nza^{2s}}dv(z)\\
&=\int_{\frac{1}{2}<|\varphi_a(z)|^2<1}\frac{2^s\nw^p\na^s\nz^{q+s}}{\nzw^{n+1+q+p}\nza^{2s}}dv(z)\\
&\leq\int_B c\frac{\na^s\nz^{s-n-1}}{\nza^{2s}}dv(z)<c,
\end{array}                                                                                                                                 \label{20}
\ee where the last two inequalities follow by theorem 1.4.10 in
\cite{Rudin} and the fact that $|1-<z,w>|\geq 1-|z|$.

And if $s\leq n$ , we choose constants $x,x',\lambda$ satisfying
$$
max\{1,\frac{n}{q+n+1}\}<x<\frac{n}{n-s}$$(when\  \ s=n, just
consider $x$ such that  $ x>max\{1,\frac{n}{q+n+1}\})$, and let $$
\lambda=q+n+1-\frac{n+1}{x}, \frac{1}{x}+\frac{1}{x'}=1,
$$
then $\lambda x>-1$ and $(q+s-\lambda)x'>-1$. Thus by
H$\ddot{o}$lder inequality \be
\begin{array}{ll}
&\int_{\frac{1}{2}<|\varphi_a(z)|^2<1}\frac{\nw^p}{\nzw^{n+1+q+p}}\nz^q\gs dv(z)\\
&\leq\int_{\frac{1}{2}<|\varphi_a(z)|^2<1}\frac{\nw^p\nz^q 2^s\na^s\nz^s}{\nzw^{n+1+q+p}\nza^{2s}}dv(z)\\
&\leq2^s\int_B\frac{\nw^p\nz^\lambda}{\nzw^{n+1+q+p}}\frac{\nz^{q+s-\lambda}\na^s}{\nza^{2s}}dv(z)\\
&\leq c\{\int_B \frac{\nw^{px}\nz^{\lambda x}}{\nzw^{(n+1+q+p)x}}dv(z)\}^\frac{1}{x}\{\int_B\frac{\nz^{(q+s-\lambda)x'}\na^{sx'}}{\nza^{2sx'}}dv(z)\}^\frac{1}{x'}\\
&\leq c .
\end{array}                                                                                                                                                                                  \label{21}
\ee

At the same time
 \be
\begin{array}{ll}
&\int_{|\varphi_a(z)|^2\leq\frac{1}{2}}\frac{\nw^p}{\nzw^{n+1+q+p}}\nz^q\gs dv(z)\\
&=\int_{|u|\leq\frac{1}{2}}\frac{\nw^p(1-|\varphi_a(u)|^2)^q\na^{n+1}}{|1-<\varphi_a(u),\omega>|^{n+1+q+p}\nua^{2n+2}}\log^s\frac{1}{|u|}dv(u)\\
&\leq c
\int_{|u|\leq\frac{1}{2}}\frac{\na^{n+1}}{(1-|\varphi_a(u)|^2)^{n+1}\nua^{2n+2}}\log^s\frac{1}{|u|}dv(u)\\
&= c\int_{|u|\leq\frac{1}{2}}\frac{1}{\nua^{n+1}}\log^s\frac{1}{|u|}dv(u)\\
&\leq c\int_{|u|\leq\frac{1}{2}}\frac{1}{\nu^{n+1}}\log^s\frac{1}{|u|}dv(u)\\
&\leq c\int_B\log^s\frac{1}{|u|}dv(u)<c .
\end{array}           \label{22}
\ee

Combine (\ref{20}),(\ref{21}),and (\ref{22}), we have
$$
\sup_{a\in B}\int_B\frac{\nw^p}{\nzw^{n+1+q+p}}\nz^q\gs dv(z)\leq c.
$$
\end{proof}

\begin{lem}There is a constant $C>0$ such that for $\forall t>-1$ and
$z\in B$ \bea
\int_B|\log\frac{1}{1-<z,w>}|^2\frac{\nw^t}{\nzw^{n+1+t}}dv(w)\leq
C\big{(}\log\frac{1}{1-|z|^2}\big{)}^{2}. \eea
\end{lem}

\begin{proof}Denote the right term as $I_t$ and let $2\lambda=t+n+1$. By
Taylor expansion
$$
 |\log\frac{1}{1-<z,w>}|^2=\sum\limits_{u,v=1}^{+\infty}\frac{<z,w>^u<w,z>^v}{uv}
$$
and
$$
\frac{1}{\nzw^{2\lambda}}=\sum\limits_{k,l=0}^{+\infty}\frac{\Gamma(\lambda+k)\Gamma(\lambda+l)}{k!l!\Gamma(\lambda)^2}<z,w>^k<w,z>^l,
$$
therefore
$$
\begin{array}{ll}
I_t&=\int_B\sum\limits_{u,v=1}^{+\infty}\sum\limits_{k,l=0}^{+\infty}\frac{\Gamma(\lambda+k)\Gamma(\lambda+l)}{uvk!l!\Gamma(\lambda)^2}<z,w>^{k+u}<w,z>^{l+v}\nw^tdv(w)\\
   &=\sum\limits_{u=1}^{+\infty}\sum\limits_{k=0}^{+\infty}\sum\limits_{l=0}^{u+k-1}\frac{\Gamma(\lambda+k)\Gamma(\lambda+l)}{u(u+k-l)k!l!\Gamma(\lambda)^2}\int_B|<z,w>|^{2(u+k)}\nw^tdv(w)
\end {array}
$$
without lost of generality, let $z=|z|e_1,$ then
$$
\begin{array}{ll}
\int_B&|<z,w>|^{2(u+k)}\nw^tdv(w)\\
&       =\int_B (|z|w_1)^{2(u+k)}\nw^tdv(w)\\
                                &=2n\int_0^1\int_{\partial B}\rho^{2n-1}|z|^{2(u+k)}|\rho\xi_1|^{2(u+k)}(1-\rho^2)^td\rho d\delta_n(\xi)\\
          &=2n|z|^{2(u+k)}\int_0^1\rho^{2(u+k+n-1)+1}(1-\rho^2)^t d\rho \int_{\partial B}|\xi_1|^{2(u+k)}d\delta(\xi)\\
          &=n|z|^{2(u+k)}\frac{\Gamma(u+k+n)\Gamma(t+1)}{\Gamma(u+k+n+t+1)}\frac{(n-1)!(u+k)!}{(u+k+n-1)!}\\
          &=\frac{\Gamma(t+1)\Gamma(u+k+1)n!}{\Gamma(2\lambda+u+k)}|z|^{2(u+k)},
\end {array}
$$
so
$$
\begin{array}{ll}
I_t&=\sum\limits_{u=1}^{+\infty}\sum\limits_{k=0}^{+\infty}\sum\limits_{l=0}^{u+k-1}\frac{\Gamma(\lambda+k)\Gamma(\lambda+l)}{u(u+k-l)k!l!\Gamma(\lambda)^2}\frac{\Gamma(t+1)\Gamma(u+k+1)n!}{\Gamma(2\lambda+u+k)}|z|^{2(u+k)}\\
   &=\sum\limits_{u=1}^{+\infty}\sum\limits_{k=0}^{+\infty}\frac{n!\Gamma(t+1)\Gamma(\lambda+k)\Gamma(u+k+1)}{uk!\Gamma(\lambda)^2\Gamma(2\lambda+u+k)}\sum\limits_{l=0}^{u+k-1}\frac{\Gamma(\lambda+l)}{(u+k-l)l!}|z|^{2(u+k)}\\
   &=\sum\limits_{u=1}^{+\infty}\sum\limits_{k=1}^{+\infty}\frac{n!\Gamma(t+1)\Gamma(\lambda+k)\Gamma(u+k+1)}{uk!\Gamma(\lambda)^2\Gamma(2\lambda+u+k)}\sum\limits_{l=0}^{u+k-1}\frac{\Gamma(\lambda+l)}{(u+k-l)l!}|z|^{2(u+k)}\\
   &\ +\sum\limits_{u=1}^{+\infty}\frac{n!\Gamma(t+1)\Gamma(u+1)}{u\Gamma(\lambda)\Gamma(2\lambda+u)}\sum\limits_{l=0}^{u-1}\frac{\Gamma(\lambda+l)}{(u-l)l!}|z|^{2u}\\
   &=I_1+I_2 ,
\end{array}
$$
by Stirling formula, there is an absolute constant $C_1$
s.t.$$\frac{\Gamma(\lambda+l)}{l!}\leq C_1 l^{\lambda-1},
\frac{\Gamma(u+k+1)}{\Gamma(2\lambda +u+k)}\leq
C_1(u+k)^{1-2\lambda}, $$$$ \frac{\Gamma(u+k+1)}{\Gamma(2\lambda
+u)}\leq C_1u^{1-2\lambda}, \frac{\Gamma(\lambda+k)}{k!}\leq C_1
k^{\lambda-1}$$ for all $l,u,k\geq 1$, then
$$
I_1\leq
C_1^3\sum\limits_{u=1}^{+\infty}\sum\limits_{k=1}^{+\infty}\frac{n!\Gamma(t+1)k^{\lambda-1}(u+k)^{1-2\lambda}}{u\Gamma(\lambda)^2}\sum\limits_{l=1}^{u+k-1}\frac{l^{\lambda-1}}{(u+k-l)}|z|^{2(u+k)}
$$
and
$$
I_2 \leq
C_1^2\sum\limits_{u=1}^{+\infty}\frac{n!\Gamma(t+1)u^{1-2\lambda}}{u\Gamma(\lambda)}\sum\limits_{l=1}^{u-1}\frac{l^{\lambda-1}}{(u-l)}|z|^{2u}
.
$$
Notice that
$$
\sum\limits_{l=1}^{M-1}\frac{l^{(\lambda-1)}}{M-l}\approx
M^{\lambda-2}\log M
$$
for any $M\geq 2$, then there is constant $C$, s.t.
$$
\begin{array}{ll}
I_1\leq  &C\sum\limits_{u=1}^{+\infty}\sum\limits_{k=1}^{+\infty} \frac{n!\Gamma(t+1)k^{\lambda-1}(u+k)^{1-2\lambda}}{\Gamma(\lambda)^2u}(u+k)^{\lambda-2}\log(u+k)|z|^{2(u+k)}\\
          &=C\sum\limits_{u=1}^{+\infty}\sum\limits_{k=1}^{+\infty}\frac{n!\Gamma(t+1)}{\Gamma(\lambda)^2}\frac{k^\lambda}{(u+k)^\lambda}\frac{\log(u+k)}{u+k}\frac{1}{uk}|z|^{2(u+k)}\\
          &\leq C\sum\limits_{u=1}^{+\infty}\sum\limits_{k=1}^{+\infty}\frac{1}{uk}|z|^{2(u+k)}=C\big{(}\log\frac{1}{1-|z|^2}\big{)}^{2}
\end{array}
$$
and
$$
\begin{array}{ll}
I_2\leq  &C\sum\limits_{u=1}^{+\infty} \frac{n!\Gamma(t+1)u^{1-2\lambda}}{\Gamma(\lambda)u}u^{\lambda-2}\log u|z|^{2u}\\
          &=C\sum\limits_{u=1}^{+\infty}\frac{n!\Gamma(t+1)}{\Gamma(\lambda)}\frac{1}{u^{\lambda+1}}\frac{\log
          u}{u}|z|^{2u},
\end{array}
$$
then it is clearly that $I_2$ can be control by
$\big{(}\log\frac{1}{1-|z|^2}\big{)}^{2}$. Therefore we are done.
\end{proof}

\begin{lem}Suppose $0<p,s<+\infty$ and $s+p>n$,then

i)\ If $s>n$, then there is a constant $C>0$, for $\forall z\in B$
$$
\begin{array}{ll}
\sup\limits_{a\in
B}\int_B\big{(}\log\frac{1}{1-|z^2|}\big{)}^{-p}\big{|}\log\frac{1}{1-<z,w>}\big{|}^p\frac{\nz^{p-n-1}}{\nzw^p}g^s(z,a)dv(z)
<C;
\end{array}
$$
ii)\ If $s\leq n$, then if choose $x$ satisfies $\max
\big{\{}1,\frac{n}{p}\big{\}}<x<\frac{n}{n-s}$  (if $n=s$, just
let $x>\max \big{\{}1,\frac{n}{p}\big{\}}$), then
$$
\begin{array}{ll}
\sup\limits_{a\in
B}\int_B\big{(}\log\frac{1}{1-|z^2|}\big{)}^{-\frac{2}{x}}\big{|}\log\frac{1}{1-<z,w>}\big{|}^\frac{2}{x}\frac{\nz^{p-n-1}}{(\nzw)^p}g^s(z,a)dv(z)<C.
\end{array}
$$
\end{lem}

\begin{proof}
 i) Notice that
$$
\begin{array}{ll}
&\big{(}\log\frac{1}{1-|z^2|}\big{)}^{-1}\big{|}\log\frac{1}{1-<z,w>}\big{|}\\
&=\big{(}\log\frac{1}{1-|z^2|}\big{)}^{-1}\big{|}\sum\limits_{k=1}^{+\infty}\frac{<z,w>^k}{k}\big{|}\\
&\leq\big{(}\log\frac{1}{1-|z^2|}\big{)}^{-1}\sum\limits_{k=1}^{+\infty}\frac{|z|^k}{k}\\
&=\big{(}\log\frac{1}{1-|z^2|}\big{)}^{-1}\log\frac{1}{1-|z^|}<1,
\end{array}
$$
and when $\frac{1}{2}<|\varphi_a(z)|^2<1$
$$
\begin{array}{ll}
\frac{\nz^{p-n-1}g^s(z,a)}{\nzw^p}\leq\frac{2^s\nz^{p+s-n-1}\na^s}{\nz^p\nza^{2s}}\leq\frac{2^{p+s}\nz^{s-n-1}\na^s}{\nza^{2s}}
\end{array}
$$
and then by(\ref{20})(\ref{22}), we see that i) holds.

ii)\ Let $\lambda=p-\frac{n+1}{x}$ and $x'=\frac{x}{1-x}$, then
$\lambda x>-1$ and $(p+s-n-1-\lambda)x'>-1$. Thus by (2.5)
$$
\begin{array}{ll}
&\int_{\frac{1}{2}<|\varphi_a(z)|^2<1}\big{(}\log\frac{1}{1-|z^2|}\big{)}^{-\frac{2}{x}}\big{|}\log\frac{1}{1-<z,w>}\big{|}^\frac{2}{x}\frac{1}{\nzw^p}\nz^{p-n-1}\times\\&\times g^s(z,a)dv(z)\\
 &\leq \int_{\frac{1}{2}<|\varphi_a(z)|^2<1}\big{(}\log\frac{1}{1-|z^2|}\big{)}^{-\frac{2}{x}}\big{|}\log\frac{1}{1-<z,w>}\big{|}^\frac{2}{x}\frac{2^s\nz^{p+s-n-1}\na^s}{\nzw^p\nza^{2s}}dv(z)\\
 &\leq 2^s \big{\{}\int_B\big{(}\log\frac{1}{1-|z^2|}\big{)}^{-2}\big{|}\log\frac{1}{1-<z,w>}\big{|}^2\frac{\nz^{\lambda x}}{\nzw^{px}}dv(z) \big{\}}^\frac{1}{x}\times\\ &\times \big{\{}\int_B\frac{\nz^{(s+p-n-1-\lambda)x'}\na^{sx'}}{\nza^{2sx'}}dv(z)\big{\}}^\frac{1}{x'}\\
&\leq C ,
\end{array}
$$
where the last inequity holds by lemma 2.6 and theorem 1.4.10 in
\cite{Rudin}.

At the same time, notice that $x>1,$  $$\begin{array}{ll}
&\int_{|\varphi_a(z)|^2\leq
\frac{1}{2}}\big{(}\log\frac{1}{1-|z^2|}\big{)}^{-\frac{2}{x}}\big{|}\log\frac{1}{1-<z,w>}\big{|}^\frac{2}{x}\frac{\nz^{p-n-1}}{\nzw^p}g^s(z,a)dv(z)\\
&\leq\int_{|\varphi_a(z)|^2\leq\frac{1}{2}}\frac{\nz^{p-n-1}}{\nzw^p}g^s(z,a)dv(z)\\
&=\int_{|u|^2\leq\frac{1}{2}}\frac{(1-|\varphi_a(u)|^2)^{p-n-1}}{|1-<\varphi_a(u),w>|^p}\frac{\na^{n+1}}{\nua^{2n+2}}\log^s\frac{1}{u}dv(u)\\
&\leq c\int_{|u|^2\leq\frac{1}{2}}\frac{1}{(1-|\varphi_a(u)|^2)^{n+1}}\frac{\na^{n+1}}{\nua^{2n+2}}\log^s\frac{1}{u}dv(u)\\
&= c\int_{|u|^2\leq\frac{1}{2}}\frac{1}{\nu^{n+1}}\log^s\frac{1}{u}dv(u)\\
&\leq c\int_B\log^s\frac{1}{u}dv(u)\leq c .
\end{array}
$$
\end{proof}

{\bf Remark}: we can also prove that when $s>n$, then
$\beta^\nqp\subset\Fs$, if we use the same method used in Lemma
2.5 and Lemma 2.7. Therefore, combine lemma 2.1, we conclude that
$\beta^\nqp=\Fs$.

\begin{lem}If $\wco:\Fs\rightarrow \nba$ is bounded, then for every
$\omega \in B $ satisfies $|\phi(\omega)|>\sqrt{\frac{2}{3}}$ and
every $u\in \mbox{\Bbbb C}^n-\{0\},$ there is a
 function $\fwu \in\Fs$  such that

$ i)$ $\exists c>0, $ independent of $\omega$ and $u$, s.t.
 $\nfwus<c ;$

$ii)$ $\{\fwu\}$ converges to zero uniformly for $u\in\mbox{\Bbbb
C}^n-\{0\}$ and z on compact subsets of $B$, when
$\phi(\omega)\rightarrow \partial B$;

$iii)$ There is a constant $c>0$, for $\forall \omega $ and $u$
\bea \|\wco \fwu\|_\nba\geq c\xzu \label{23} \eea holds. \end{lem}

\begin{proof} First we suppose $\phi(\omega)=\rw e_1,$ where $\rw=|\pw|$,
$e_1$ is the vector\\$(1, 0,0,\cdots,0)$.

If $|\pw|>\sqrt{2/3}$ and $\sqrt{\npw}|\Jpu|\leq|<\pw,\Jpu>|$.
Take \be \fwu(z)=\frac{(z_1-\rw)(\nrw)}{(1-\rw z_1)^{\nqp+1}}.
\label{24} \ee

Then
$$
\frac{\partial\fwu(z)}{\partial z_1}=\frac{\nrw}{(1-\rw
z_1)^{\nqp+1}}\Big{(}1+\nqp\frac{(z_1-\rw)\rw}{1-\rw z_1}\Big{)}
$$
and
$$
\frac{\partial\fwu(z)}{\partial z_k}=0, \  \ k=2,\cdots,n.
$$
Therefore \be
\begin{array}{ll}
|\nabla\fwu(z)|&=\frac{\nrw}{|1-\rw z_1|^{\nqp+1}}\Big{|}1+\nqp\frac{(z_1-\rw)\rw}{1-\rw z_1}\Big{|}\\
&\leq (1+\nqp)\frac{\nrw}{|1-\rw z_1|^{\nqp+1}}.
\end{array}                                                                                                                      \label{25}
\ee

By Lemma 2.5, $\fwu\in\Fs$, and exists $c>0$ independent of
$\omega$ and $u$, s.t. $\|\fwu\|\leq c. $ or $i)$ holds.

On the other hand, by (\ref{2}) \be
\begin{array}{ll}
&\|\wco\fwu\|_\nba\geq \nw^\alpha|\nabla(\wco\fwu)(\omega)|\\
                   &=\nw^\alpha|\psi(\omega)\nabla(\fwu\circ\phi)(\omega)|,\end{array}\ee
notice that
$$
<\phi(\omega),\Jpu>=\rw e_1\Jpu,$$ and use lemma 2.4, we see that
$$\begin{array}{ll}
                 &\|\wco\fwu\|_\nba\geq c\nw^\alpha|\psi(\omega)|\Big{|}\frac{\nabla\fwu(\pw)\Jpu}{\sqrt{\nw|u|^2+|<\omega,u>|^2}}\Big{|}\\
                   &=c\nw^\alpha|\psi(\omega)|\frac{\nrw}{(\nrw)^{\nqp+1}}\Big{|}\frac{e_1\Jpu}{\sqrt{\nw|u|^2+|<\omega,u>|^2}}\Big{|}\\
                   &=c|\psi(\omega)|\frac{\nw^\alpha}{\rw(\nrw)^{\nqp}}\Ba\frac{<\pw,\Jpu>}{\wu}\Ba,
                   \end{array}$$
                   therefore, by our assumption, we get
                   $$\begin{array}{ll}
                   &\|\wco\fwu\|_\nba\geq c|\psi(\omega)|\frac{\nw^\alpha}{(\nrw)^{\nqp}}\frac{\sqrt{|<\pw,\Jpu>|^2+|<\pw,\Jpu>}|}{\wu}\\
                   &= c
                   \xzu .
\end{array}                                                                                                                                                          \label{26}
$$

If $\sqrt{\npw}|\Jpu|\leq|<\pw,\Jpu>|$. Write $J_\pw
u=\\(\xi_1,\cdots,\xi_n)^T.$ For $j=2,\cdots,n,$ let
$\theta_j=\arg\xi_j$ and $a_j=e^{-i\theta_j}$ (if $\xi_j=0$, let
$a_j=0$).

Let \be \fwu(z)=\frac{(a_2z_2+\cdots+a_nz_n)(\nrw)^{3/2}}{(1-\rw
z_1)^{\nqp+1}}, \label{27} \ee then
$$
\frac{\partial\fwu(z)}{\partial
z_1}=\frac{(\nqp+1)\rw(\nrw)^{3/2}}{(1-\rw
z_1)^{\nqp+2}}(a_2z_2+\cdots+a_nz_n)
$$
and
$$
\frac{\partial\fwu(z)}{\partial z_k}=\frac{a_k(\nrw)^{3/2}}{(1-\rw
z_1)^{\nqp+1}}, \  \ k=2,\cdots,n.
$$

Therefore \be
\begin{array}{ll}
|\nabla\fwu(z)|&=\sqrt{|\frac{\partial\fwu(z)}{\partial z_1}|^2+|\frac{\partial\fwu(z)}{\partial z_2}|^2+\cdots+|\frac{\partial\fwu(z)}{\partial z_1}|^2}\\
               &=\sqrt{\frac{(\nqp+1)^2\rw^2(\nrw)^3|a_2z_2+\cdots+a_nz_n|^2}{|1-\rw z_1|^{2(\nqp+2)}}+\frac{(n-1)(\nrw)^3}{|1-\rw z_1|^{2(\nqp+1)}}}\\
               &\leq\sqrt{\frac{(n-1)(\nqp+1)^2\rw^2(\nrw)^3(|z_2|^2+\cdots+|z_n|^2)}{|1-\rw z_1|^{2(\nqp+2)}}+\frac{(n-1)(\nrw)^3}{|1-\rw z_1|^{2(\nqp+1)}}}\\
               &\leq\frac{\sqrt{(n-1)}(\nrw)^{3/2}}{|1-\rw z_1|^{\nqp+1}}\sqrt{\frac{(\nqp+1)^2\rw^2(1-|z_1|^2)}{|1-\rw z_1|^2}+1}\\
               &\leq\frac{\sqrt{(n-1)}(\nrw)}{|1-\rw z_1|^{\nqp+1}}|\nrw+\rw^2(\nqp+1)^2)|^\frac{1}{2}\\
               &\leq c\frac{\nrw}{|1-\rw z_1|^{\nqp+1}},
               \end{array}                                                                                                                                                             \label{28}
\ee so, by lemma 2.5, we know that $\fwu\in\Fs,$ and $\fwu$
satisfies $i)$. Next we will extimate $\|\wco\fwu\|_\alpha$.\\
Since $|\pw|>\sqrt{2/3}$ and $$\sqrt{\npw}|\Jpu|>|<\pw,\Jpu>|,$$
we can get$$\sqrt{3(\npw)(|\xi_2|^2+\cdots+|\xi_n|^2)}>|\xi_1|.$$
Then
$$|\xi_2|^2+\cdots+|\xi_n|^2\geq\frac{1}{2}(|\xi_1|^2+\cdots+|\xi_n|^2).$$
Therefore \be
\begin{array}{ll}
&\|\wco\fwu\|_\alpha\geq c\nw^\alpha|\psi(\omega)|\Big{|}\frac{\nabla\fwu(\pw)\Jpu}{\sqrt{\nw|u|^2+|<\omega,u>|^2}}\Big{|}\\
                   &=c\nw^\alpha|\psi(\omega)|\frac{\nw^\alpha(\nrw)^\frac{1}{2}(|\xi_2|+\cdots+|\xi_n|)}{(\nrw)^{\nqp}\wu}\\
                   &\geq c\nw^\alpha|\psi(\omega)|\frac{\nw^\alpha(\nrw)^\frac{1}{2}\sqrt{|\xi_2|^2+\cdots+|\xi_n|^2}}{(\nrw)^{\nqp}\wu}\\
                   &\geq c\nw^\alpha|\psi(\omega)|\frac{\nw^\alpha(\nrw)^\frac{1}{2}\sqrt{|\xi_1|^2+|\xi_2|^2+\cdots+|\xi_n|^2}}{(\nrw)^{\nqp}\wu}\\
                   &\geq c\nw^\alpha|\psi(\omega)|\frac{\nw^\alpha(\nrw)^\frac{1}{2}|\Jpu|}{(\nrw)^{\nqp}\wu}\\
                   &\geq c\xzu.
\end{array}
\ee \label{29}

By our discussion above, for every $\omega \in B$ satisfies
$|\pw|>\sqrt{2/3}$ and $u\in\mbox{\Bbbb C}^n-\{0\}$, there is a
$\fwu\in \Fs$ satisfies $i)$ and $iii)$. And that $\fwu$ has the
property $iii)$ is clearly, since it is very easy to verify that
for every
 $\omega$ and $u$
 $$
 |\fwu(z)|\leq c\frac{\nrw}{\nz^{\nqp+1}}
 $$
holds.

In general situation, or if $\pw\neq|\pw|e_1,$ we use the unitary
transformation $\Uw$ to make $\pw=\rw e_1\Uw$, where $\rw=|\pw|$.
Then $\gwu=\fwu\circ U_\omega^{-1}$ is the desired function.

In fact, by $\nabla \gwu(z)=\nabla(\fwu\cdot
U_\omega^{-1})(z)=(\nabla\fwu)(zU_\omega^{-1})(U_\omega^{-1})^T$
and
   \\$|zU_\omega^{-1}|=|z|$, we have
$$
\begin{array}{ll}
&\int_B|\nabla\gwu(z)|^p\nz^q\gs dv(z)\\
&=\int_B|(\nabla\fwu)(zU_\omega^{-1})(U_\omega^{-1})^T|^p\nz^q\gs dv(z)\\
&=\int_B|\nabla\fwu(z)|^p\nz^q\gs dv(z),
\end{array}
$$
where the last equation we use the linear coordinate translation
$z=zU_\omega^{-1}$ and the fact that $\Fs$ is invariant under
$m\ddot{o}bious$ translation. So $$\nfwus=\ngwus.$$

Then we can prove the same result by the same way,  we omit the
details here.
\end{proof}

\begin{lem}Let $0\leq p <1,$ then for every $f\in\beta^p$ and $\forall
z,w\in B$ satisfy $z=rw$ for some real number $r$, \be
|f(z)-f(w)|\leq\frac{2\|f\|_p}{1-p}|z-w|^{1-p} \label{30} \ee
Especially  when $0<\nqp<1,$ there is a constant $c>0$ such that \be
|f(z)-f(w)|\leq c\frac{\|f\|_\Fs}{1-p}|z-w|^{1-p} \label{31} \ee
 for $\forall f\in\Fs$ and $z,w$ as above.
\end{lem}

\begin{proof}Denote $F(t)=f(tz+(1-t)w)$. Then
$$
\begin{array}{ll}
f(z)-f(w)&=F(1)-F(0)\\
         &=\int^1_0 F'(t)dt\\
         &=\int^1_0\sum^n_{k=1}\frac{\partial f(tz+(1-t)w)}{\partial \xi_k}(z_k-w_k)dt\\
         &=\int^1_0<\nabla f(tz+(1-t)w),\overline{z-w}>dt,
\end{array}
$$
so
$$
\begin{array}{ll}
|f(z)-f(w)|&=|\int^1_0<\nabla f(tz+(1-t)w),\overline{z-w}>dt|\\
           &\leq \int^1_0|\nabla f(tz+(1-t)w)||z-w|dt.
\end{array}
$$
Since $f\in \beta^p,$ it's easy to show that
$$
|\nabla
f(z)|\leq\frac{\|f\|_{\beta^p}}{\nz^p}\leq\frac{\|f\|_{\beta^p}}{(1-|z|)^p},
\  \ \forall z\in B.
$$
Note that if two numbers $a,b$, satisfy $|a+b|<1,|a|<1$ and $|b|<1$,
\\ then $|a+b|+|a|+|b|<2$ and so $1-|a+b|>|(1-|a|-|b|)|$. And notice that $z=rw$ for some real number $r$,
thus
$$
\begin{array}{ll}
|f(z)-f(w)|&\leq \int^1_0\frac{\|f\|_{\beta^p}|z-w|}{(1-|w+t(z-w)|)^p}dt\\
           &\leq\int^1_0\frac{\|f\|_{\beta^p}|z-w|}{|1-|w|-t|z-w||^p}dt.
\end{array}
$$
If $1-|w|-|z-w|\geq 0$, then \be\begin{array}{ll}
|f(z)-f(w)|&\leq\int^1_0\frac{\|f\|_{\beta^p}|z-w|}{(1-|w|-t|z-w|)^p}dt\\
           &\leq\frac{\|f\|_{\beta^p}}{1-p}(1-|w|-t|z-w|)^{(1-p)}\mid^0_1\\
           &=\frac{\|f\|_{\beta^p}}{1-p}\{(1-|w|)^{1-p}-(1-|w|-|z-w|)^{1-p}\}.\label{2.19}
\end{array}\ee
Write $\varphi(x)=x^{1-p}-b^{1-p}-(x-b)^{1-p}, \  \ x\in(b,a)$, then
$$
\varphi
'(x)=(1-p)x^{-p}-(1-p)(x-b)^{-p}=(1-p)(\frac{1}{x^p}-\frac{1}{(x-b)^p})<0,
$$
since $b<x<a$. Therefore $\varphi(a)<\varphi(b)=0$, or
$a^{1-p}-b^{1-p}-(a-b)^{1-p}<0$, if $b<a$. If we let $a=1-|w|,
b=|z-w|, $ then
$$
(1-|w|)^{1-p}<|z-w|^{1-p}+(1-|w|-|z-w|)^{1-p}.
$$
Combine (\ref{2.19}), we have
$$
|f(z)-f(w)|\leq \frac{\|f\|_p}{1-p}|z-w|^{1-p} .
$$

If $1-|w|-|z-w|<0$, then there is a $\xi\in (0,1),$ s.t.
$$
1-|w|-\xi|z-w|=0.
$$
so
$$
\begin{array}{ll}
|f(z)-f(w)|&\leq\int^\xi_0\frac{\|f\|_{\beta^p}|z-w|}{(1-|w|-t|z-w|)^p}dt+\int^1_\xi\frac{\|f\|_{\beta^p}|z-w|}{(|w|+t|z-w|-1)^p}dt\\
           &\leq\frac{\|f\|_{\beta^p}}{1-p}\{(1-|w|-t|z-w|)^{(1-p)}\mid^0_\xi+(|w|+t|z-w|)^{(1-p)}\mid^1_\xi\}\\
           &=\frac{\|f\|_{\beta^p}}{1-p}\{(1-|w|)^{1-p}+(|w||z-w|-1)^{1-p}\}\\
          &\leq \frac{2}{1-p}\|f\|_{\beta^p}|z-w|^{1-p}.
\end{array}
$$

Therefore for all $z,w\in B,$
$$
|f(z)-f(w)|\leq\frac{2\|f\|_{\beta^p}}{1-p}|z-w|^{1-p}.
$$

The second result of this lemma follows quickly by lemma 2.1.
\end{proof}

\section{The proof of theorem 1.1}

We'll give the sufficiency first. Suppose that (\ref{5}) and
(\ref{6}) hold.
\be
\begin{array}{ll}
&\nw^\alpha|\nabla\wco f(\omega)|\\
 &=\nw^\alpha|\psi(\omega)\nabla(f\circ\phi)(\omega)+f(\pw)\nabla\psi(\omega)|\\
 &\leq\nw^\alpha|\psi(\omega)\nabla(f\circ\phi)(\omega)|+\nw^\alpha|f(\pw)\nabla\psi(\omega)|\\
&=T_1+T_2 .\end{array}\ee By lemma 2.4 and (1.2),
$$
\begin{array}{ll}
T_1&\leq c\sup\limits_{u\in\mbox{\Bbbb C}^n-\{0\}}\frac{\nw^\alpha|\phi(\omega)||\nabla f(\nw)\Jpu|}{\wu}\\
&=c\sup\limits_{u\in\mbox{\Bbbb C}^n-\{0\}}\frac{\nwa|\psi(\omega)|}{\npwq}\frac{\Jpwu}{\wu}\times\\
   &\  \ \frac{\npwq\|\nabla f(\pw)\cdot\Jpu|}{\Jpwu}+\\
   &\leq c\sup\limits_{u\in\mbox{\Bbbb C}^n-\{0\}}\xzu\times \\
   &\  \ \|f\|_{\beta^\nqp},
\end{array}
$$
by lemma 2.1 and (\ref{5}), we have
$$
T_1\leq c\|f\|_{\beta^\nqp}.
$$
On the other hand, by lemma 2.2 and (\ref{6}),
$$\begin{array}{ll}
T_2&\leq c\nwa G_\nqp(\pw)\| f\|_{F(p,q,s)}|\nabla \psi(\omega)|\\
&\leq c\|f\|_{\beta^\nqp},
\end{array}$$
notice that \be |\wco
f(0)|=|\psi(0)||f(\phi(0))|\leq|f(\phi(0))|\leq G_\nqp(\phi(0))\|
f\|_{F(p,q,s)},  \ee therefore, we have
$$
\|\wco f\|_\nba\leq    c\| f\|_{F(p,q,s)}.
$$
Since the constant $c$ is independent of $f$,
$\wco:\Fs\rightarrow\nba$ is bounded.

One the other hand, suppose $\wco$ is bounded, with
$$
\|\wco f\|_\nba\leq    c\| f\|_{F(p,q,s)},
$$
for all $f \in\Fs.$ It's very easy to show that the functions
$f_l(z)=z_l,\  \ l=1,\cdots, n$ and $f(z)=1$ are in $\Fs$.
Therefore, $\wco f_l$ and $\wco f$ must in $\nba$ too, in an other
word,  $\psi\phi_l, l=1,\cdots,n$ and $\psi$ lies in $\nba$.

(1) If $|\pw|^2\leq 2/3$, then \be
\begin{array}{ll}
&\xzu\\
&\leq c|\psi(\omega)|\nwa\frac{|\Jpu|}{\wu}\\
&\leq c|\psi(\omega)|\nwa\{\sum^n_{l=1}\frac{|\nabla \phi_l(\omega)u|}{\wu}\}\\
&\leq c\nwa\{\sum^n_{l=1}|\psi(\omega)\nabla \phi_l(\omega)|\}\\
&= c\nwa\{\sum^n_{l=1}[|\nabla(\psi\cdot\phi_l)(\omega)-\phi_l(\omega)\cdot\nabla\psi(\omega)|]\}\\
&\leq c\nwa\{\sum^n_{l=1}[|\nabla(\psi\cdot\phi_l)(\omega)|+|\phi_l(\omega)\cdot\nabla\psi(\omega)|]\}\\
&\leq c\sum^n_{l=1}\|\psi\phi_l\|_{\nba}+\|\psi\|_{\nba}\\
&\leq c ,
\end{array}                                                                                                                     \label{33}
\ee where the last inequality follows by the fact that
$\psi\cdot\phi_l(l=1,\cdots,n)$ and $\psi(z)$ lies in $\nba$.

(2) If $|\pw|^2> 2/3$, then by lemma 2.8, for every $\omega \in B
$ satisfies $|\phi(\omega)|^2>2/3$ and every $u\in \mbox{\Bbbb
C}^n-\{0\},$ exists a
 function $\fwu \in\Fs$  such that $ \nfwus<c $, and
 $$
 \begin{array}{ll}
&\|\wco \fwu\|_\nba \\
&\geq c\xzu,
\end{array}
$$
but $\|\wco\fwu\|_\nba\leq\|\wco\|\|\fwu\|_\Fs,$ we get \bea&
\xzu\\&\leq c. \label{34} \eea

 By (\ref{33}) and (\ref{34}), we know
that (\ref{5}) holds.

  If $0<\nqp<1$, it's obviously that (\ref{6}) holds, since $\psi$
lies in $\nba$. By the previous discussion, we
 can suppose that $\pw=\rw e_1$.

   If $\nqp>1,$ take
   $$
\begin{array}{ll}
   \fw(z)=\frac{\nrw}{(1-\rw z_1)^{\nqp}}                                                                                                \label{36}
   \end{array}$$
   then
   $$
\frac{\partial\fw(z)}{\partial z_1}=\nqp\frac{\rw(\nrw)}{(1-\rw
z_1)^{\nqp+1}}
$$
and
$$
\frac{\partial\fw(z)}{\partial z_k}=0, \  \ k=2,\cdots,n.
$$
so, by lemma 2.5, $\fw\in \Fs$, and
$$
\|\fw\|_{\Fs}\leq c,
$$
for every $\omega$ in $B$. Therefore \be
\begin{array}{ll}
&\nwa  |\nabla\psi(\omega)|G_{\nqp}(\pw)=\frac{\nwa|\nabla\psi(\omega)|}{(\nrw)^{\nqp-1}}\\
            &=\nwa|\nabla\psi(\omega)||\fw(\pw)|\\
            &\leq\nwa\{|\nabla(\psi\cdot \fw\circ\phi)(\omega)|+|\psi(\omega)\nabla(\fw\circ\phi)(\omega)|\}\\
            &\leq \|\wco
            \fw\|_\nba+\nwa|\psi(\omega)\nabla(\fw\circ\phi)(\omega)|,
\end{array}\ee
since $\wco$ is bounded, we have $$\|\wco \fw\|_\nba\leq
c\|\fw\|_{\Fs}\leq c.$$ And by lemma 2.4, (\ref{5}) and our
assumption
$$
\begin{array}{ll}
&\nwa|\psi(\omega)\nabla(\fw\circ\phi)(\omega)|\\
            &\leq c\sup\limits_{u\in\mbox{\Bbbb C}^n\setminus\{0\}}\frac{\nwa|\psi(\omega)||\nabla\fw(\pw)\circ\Jpu|}{\wu}\\
            &\leq c\nqp\sup\limits_{u\in\mbox{\Bbbb C}^n\setminus\{0\}}|\psi(\omega)|\frac{\nwa}{(\nrw)^\nqp}\frac{\sqrt{|<\pw,\Jpu>|^2}}{\wu}\\
            &\leq c\Xzu\\
            &\leq c.
\end{array}                                                                                                                                  \label{37}
$$

If $\nqp=1$ and $s>n$ for every $\omega\in B$,take
$$
\fw(z)=\Big{(}\log\frac{1}{\nrw}\Big{)}^{-1}\Big{(}\log\frac{1}{1-\rw
z_1}\Big{)}^2
$$
then
$$
\frac{\partial\fw(z)}{\partial
z_1}=\Big{(}\log\frac{1}{\nrw}\Big{)}^{-1}\log\frac{1}{1-\rw
z_1}\frac{\rw}{1-\rw z_1},
$$
and
$$
\frac{\partial\fw(z)}{\partial z_l}=0, \  \ l=2,\cdots,n.
$$
Therefore
$$
|\nabla\fw(z)|=\Big{(}\log\frac{1}{\nrw}\Big{)}^{-1}\big{|}\log\frac{1}{1-\rw
z_1}\big{|}\frac{\rw}{|1-\rw z_1|},.
$$
so, by lemma 2.7, $\fw\in \Fs$, and there is a constant c, s.t.
$$
\|\fw\|_{\Fs}\leq c,
$$
for every $\omega\in B$. By the same calculation the above
discussion,
$$
\begin{array}{ll}
&\nwa|\nabla\psi(\omega)|\log\frac{1}{1-\rw^2}\\
&=\nwa|\nabla\psi(\omega)||\fw(\pw)|\\
            &\leq \|\wco \fw\|_\nba+\nwa|\psi(\omega)\nabla(\fw\circ\phi)(\omega)|\\
            &\leq c\|\fw\|_{\Fs}+c\sup\limits_{u\in\mbox{\Bbbb C}^n-\{0\}}\frac{\nwa|\psi(\omega)||\nabla\fw(\pw)\circ\Jpu|}{\wu}\\
            &= c+c\nqp\sup\limits_{u\in\mbox{\Bbbb C}^n-\{0\}}|\psi(\omega)|\frac{\nwa(1-\rw^{k-1})}{(\nrw)(1-\rw^{k+1})}\frac{\sqrt{|<\pw,\Jpu>|^2}}{\wu}\\
            &\leq c+c\sup\limits_{u\in\mbox{\Bbbb C}^n-\{0\}}|\psi(\omega)|\frac{\nwa}{\nrw}\frac{\sqrt{|<\pw,\Jpu>|^2}}{\wu}\\
            &\leq c+c\sup\limits_{u\in\mbox{\Bbbb C}^n-\{0\}}|\psi(\omega)|\frac{\nwa}{\nrw}\frac{\Jpwu}{\wu}\\
            &\leq c.
\end{array}
$$
Notice that $\psi\in\nba$, so we have
$$
\nwa|\nabla\psi(\omega)|\log\frac{2}{\nrw}<\infty
$$
or
$$
\nwa|\nabla\psi(\omega)|G_1(\pw)<\infty
$$
\\
when $s\leq n$, just take
$$
\fw(z)=\Big{(}\log\frac{1}{\nrw}\Big{)}^{-\frac{2}{px}}\Big{(}\log\frac{1}{1-\rw
z_1}\Big{)}^{1+\frac{2}{px}} ,
$$
where $x$ is defined in lemma 2.7, then use lemma 2.7 and by the
same proof as the case $s>n$, we can proof the same result.

Now the proof of Theorem 1 is completed.

\section{the Proof of Theorem 1.2}

If $\wco$ is bounded and (\ref{8}) holds.  ${f_j}$ is a sequence in
$\Fs$ which converges to
 zero uniformly on compact subsets of $B$ and $$\|f_j\|_{\Fs}\leq
 1$$
 for every $j\in\mbox{\Bbbb N}$.  $\forall \varepsilon >0,$ suppose $\delta=\varepsilon^{\frac{1}{1-\nqp}}$, there is a
 $k_0\in\mbox{\Bbbb N},$  for $k>k_0,$ and $\forall w\in \overline{B_{1-\delta}}=\{z:\  \ |z|\leq1-\delta\}$ have
 $$
 |f_k(w)|<\varepsilon.
 $$
 One the other hand, by lemma 2.9, if $w=rz$ for some real number $r$, $|z-w|<2\delta$ and
 $1-|w|<\delta$,
\be
 |f_k(z)-f_k(w)|\leq c\frac{\|f_k\|_\Fs}{1-p}|z-w|^{1-\nqp}\leq c\varepsilon
\ee
 holds for every $k\in\mbox{\Bbbb N}$. For every $w\in B\setminus\overline{B_{1-\delta}} $ or in a other word, $1-\delta<|w|<1$,                    \label{38}
 we can always find a point  $z_w\in B_{1-\delta}$ s.t. $w=rz$ for some real number $r$ and $|z-w|<2\delta$, then if $k>k_0$
 $$
 \begin{array}{ll}
 |f_k(w)|&\leq |f_k(w)-f_k(z_w)|+|f_k(z_w)|\leq c\varepsilon +\varepsilon<c\varepsilon.
 \end{array}
 $$

 Then, $\forall\varepsilon>0$, there is a $k_0>0,$ s.t. for all $z\in B $ and $k>k_0$,
 $|f_k(z)|<\varepsilon$ holds. So, the sequence $\{f_j\}$ converges to zero uniformly on $B$.

 Since (\ref{8}) holds, therefore, there is a $\delta>0$, such that if $1-|\pw|>\delta$
 \bea
 \sup\limits_{u\in\mbox{\Bbbb C}^n\setminus\{0\}}\xzu<\varepsilon.                                                  \label{39}
 \eea

So, when  $1-|\pw|>\delta$, \be
\begin{array}{ll}
&\nwa|\nabla\wco f_k(\omega)|\\
&\leq\nwa|\nabla\psi(\omega)f_k(\pw)|+\nwa|\psi(\omega)\nabla (f_k\circ\phi)(\omega)|\\
&\leq \|\psi\|_\nba |f_k(\pw)|+c\sup\limits_{u\in\mbox{\Bbbb C}^-\{0\}}\nwa|\psi(\omega)|\frac{|\nabla f_k(\pw)\cdot \Jpu|}{\wu}\\
&\leq c\varepsilon+c\sup\limits_{u\in\mbox{\Bbbb C}^n-\{0\}}\frac{\nwa|\psi(\omega)|}{\npwq}\frac{\Jpwu}{\wu}\times\\
                                  &\frac{\npwq\|\nabla f(\pw)\cdot\Jpu|}{\Jpwu}\\
                                  &\leq c\varepsilon+ \\
                                  &c\sup\limits_{u\in\mbox{\Bbbb C}^n-\{0\}}\xzu\\
                                  & \  \ \times\|f\|_{\Fs}\\
                                  &\leq c\varepsilon .
\end{array}         \label{40}
\ee

When $|\pw|\leq\delta,$ let $E=\{\omega:|\omega\leq\delta|\}$, then
E is compact and $\{\omega\in B:\pw\leq\delta\}\subset E$. Therefore
\be
\begin{array}{ll}
&\nwa|\nabla\wco f_k(\omega)|\\
&\leq\nwa|\nabla\psi(\omega)f_k(\pw)|+\nwa|\psi(\omega)\nabla (f_k\circ\phi)(\omega)|\\
&\leq\nwa|\nabla\psi(\omega)f_k(\pw)|+\nwa|\psi(\omega)\nabla f_k(\pw)||J_{\phi}(\omega)|\\
&\leq \psi\|_\nba |f_k(\pw)|+\sum_{l=1}^n(\|\psi\phi_l\|_\nba+\|\psi\|_\nba)|\nabla f_k(\pw)|\\
&<c\varepsilon.
\end{array}\label{41}
\ee where the last inequality follows by the fact that the sequence
$\{\nabla f_k\}$ converges to zero uniformity on $E$.

Therefore $\nwa|\nabla\wco f_k(\omega)|$ converges to zero uniformly
on $B$.

Notice that $\{f_k(0)\}$ converges to zero, and so we have
$$
\|\wco f_k\|_\nba\rightarrow 0
$$
as $|\pw|\rightarrow 0$.

On the other hand, if $\wco $ is compact, then clearly that $\wco $
is bounded and by lemma 2.8 and lemma 2.3, it is easy to know that
(\ref{8}) holds. And so we are done.

\section{The proof of theorem 1.3}

Just as the proof of theorem 2, if $\wco$ is bounded and
(\ref{9})(\ref{10}) hold,  ${f_j}$ is a sequence in $\Fs$ which
converges to
 zero uniformly on compact subsets of $B$ and $\|f_j\|_{\Fs}\leq 1$
 for every $j\in\mbox{\Bbbb N}$, then by the same discussion
as lemma 2.8, we have \be
\begin{array}{ll}
&\nw^\alpha|\nabla\wco f_j(\omega)|\\
&=\nw^\alpha|\psi(\omega)\nabla(f_j\circ\phi)(\omega)+f_j(\pw)\nabla\psi(\omega)|\\
                                  &\leq c\sup\limits_{u\in\mbox{\Bbbb C}^n-\{0\}}\xzu \times
                                  \\ &\times \|f_j\|_{\Fs}+ c\nwa G_\nqp(\pw)|\nabla \psi(\omega)|\| f_j\|_{F(p,q,s)}\\
                                  &\leq c\sup\limits_{u\in\mbox{\Bbbb C}^n-\{0\}}\xzu + \\
                                  &c\nwa\nabla \psi(\omega)|
                                  G_\nqp(\pw)|.
\end{array}    \label{42}
\ee Therefore, $\forall\varepsilon>0,$ there is a $\delta>0,$ for
$\forall j\in\mbox{\Bbbb N}$ and $1-|\pw|<\delta$, \be
\nw^\alpha|\nabla\wco f_j(\omega)|\leq c\varepsilon\label{43} \ee
follows.

And when $1-|\pw|\leq\delta$, by the same discussion as the proof
of Theorem 2, we know that there is a $k_0\in\mbox{\Bbbb N},$ s.t.
when $j>k_0$ have \be \nw^\alpha|\nabla\wco f_j(\omega)|\leq
c\varepsilon\label{44}. \ee

Combine (\ref{43}) and (\ref{44}), and notice that $\{f_k(0)\}$
converges to zero, we have
$$
\|\wco f_k\|_\nba\rightarrow 0,
$$
as $|\pw|\rightarrow 1$, therefore by lemma 3, $\wco$ is compact.

If  $\wco$ is compact, then $\wco$ is bounded, and by lemma 2.8 and
lemma 2.3, it's clearly that (\ref{9}) hold.

If $\nqp=1, s>n$, and $|\rw|>\frac{2}{3}$, let \be
\fw=\Big{(}\log\frac{1}{\nrw}\Big{)}^{-1}\Big{(}\log\frac{1}{1-\rw
z_1}\Big{)}^2,\label{45} \ee by the proof of theorem 1,
$\fw\in\Fs,$ and $\exists C>0$, such that $\|\fw\|_ {\Fs}<C$. And
it is clearly that {$\fw$} converges to zero uniformly  on compact
subsets of $B$ by the same discussion as the proof of theorem 1.1
and use lemma 2.1,we have \be
\begin{array}{ll}
&\nwa|\nabla\psi(\omega)||\fw(\pw)|\\
            &=\nwa|\nabla\psi(\omega)|\log\frac{1}{1-\rw^2}\\
                       &\leq c\|\wco \fw\|_{\nba}\\
                       &+c\sup_{u\in\mbox{\Bbbb
                       C}^n-\{0\}}|\psi(\omega)|\frac{\nwa}{\nrw}\frac{\Jpwu}{\wu},
          \end{array}\label{46}
\ee but by lemma 2.3 and (\ref{9}), we get
$$
\nwa|\nabla\psi(\omega)|\log\frac{1}{1-\rw^2}\rightarrow 0,
$$
as $|\pw|\to 1.$ Notice that (\ref{6}) holds, so
$$\nwa|\nabla\psi(\omega)|\to 0,\hspace{4mm}|\pw|\to 1.$$ And then
$\nwa|\nabla\psi(\omega)|\log\frac{2}{1-\rw^2}\to 0 $ (  $|\pw|\to
1$.)

When $\nqp=1, s<n$, just take
$$
\fw(z)=\Big{(}\log\frac{1}{\nrw}\Big{)}^{-\frac{2}{px}}\Big{(}\log\frac{1}{1-\rw
z_1}\Big{)}^{1+\frac{2}{px}},
$$
where $x$ is the one used in lemma2.6. The same discussion as the
case $s>n$ gives the same result,and we omit it here.

If $\nqp>1$, just let
$$
\fw=\frac{1-\rw^2}{(1-\rw z_1)^{\nqp+1}}
$$
and use the same method as the situation of  $\nqp<1$, we can also
prove that (\ref{10}) holds. So, the proof of theorem 1.3 is
completed.

\end{document}